\documentclass{amsproc}
\usepackage{enumerate}
\usepackage{stmaryrd}
\usepackage{dsfont}

\usepackage{tikz} 
\usetikzlibrary{matrix,arrows}

\newcounter{thm}
\newtheorem{definition}[thm]{Definition}
\newtheorem{conjecture}[thm]{Conjecture}
\theoremstyle{definition}
\newtheorem{example}[thm]{Example}

\newcommand\sheaf{\mathcal}
\newcommand\CC{\mathbb C}
\newcommand\PP{\mathbb P}
\newcommand\ZZ{\mathbb Z}
\newcommand\trip{\ensuremath{ \sheaf O_1 \sheaf O_2 \sheaf O_3}}
\newcommand\modcpt{{\ensuremath{\overline{M(X,\beta)}}}}
\newcommand{\takes}[2]{\!:\!#1 \rightarrow #2}

\newcommand{\wt}{\widetilde}

\usepackage[xetex, bookmarks, colorlinks, breaklinks,pdfauthor={Josh Guffin}]{hyperref}  
\hypersetup{linkcolor=blue,citecolor=blue,filecolor=black,urlcolor=blue} 

\begin{document}
\bibliographystyle{halpha}

\title{Quantum Sheaf Cohomology, a pr\'ecis}
\author{Josh Guffin}
\address{Department of Mathematics, University of Pennsylvania,
Philadelphia, PA 19104}
\email{guffin@math.upenn.edu}

\date{5 January, 2011}

\begin{abstract}
  We present a brief introduction to quantum sheaf cohomology, a 
  generalization of quantum cohomology based on the physics of the (0,2)
  nonlinear sigma model.
\end{abstract}

\maketitle
\thispagestyle{empty}
\newcommand{\ext}[1]{\ensuremath{{\textstyle\bigwedge}^{\! #1 }}}

This paper is based on a talk given on 13 December, 2010 during the
\emph{Second Latin Congress on Symmetries in Geometry and Physics} at the
Universidade Federal do Paran\'a in Curitiba, Brazil.  Throughout, we will
consider $X$ to be a K\"ahler manifold of complex dimension $n$.  In
addition, we will consider $\sheaf E \to X$ to be a complex Hermitian
vector bundle of rank $k$ satisfying
\begin{enumerate}[(i)]
  \item $c_2(\sheaf E) = c_2(T_X)$, \label{cond:c2}
  \item $\det \sheaf E^\vee \cong \omega_X$\label{cond:det}.
\end{enumerate}
As these conditions imply the usual Green-Schwarz anomaly cancellation
conditions, we will call such a bundle \emph{omalous}\footnotemark.  One
may consider \eqref{cond:det} to be an analogue of the usual condition for
existence of the $B$-model.  A bundle satisfying these conditions may be
obtained by, for example, selecting a deformation of the tangent bundle
when $X$ is a projective toric variety.

\footnotetext{That is, \emph{not anomalous} -- this delightful terminology is
due to Ron Donagi.}

\section*{Quantum Cohomology}

\subsection*{Ordinary cohomology}
We now give some elementary facts about the cohomology of $X$, stated in a
way that will facilitate our point of view on quantum sheaf cohomology.
Since $X$ is K\"ahler, there is a Hodge decomposition on $H^\bullet(X, \CC)$, 
\[
H^\bullet(X, \CC) \cong \bigoplus_{p,q} H^p(X, \ext q T_X^\vee).
\]
By a slight abuse of language, we will refer to elements of
the sheaf cohomology groups $H^p(X, \ext q T_X^\vee)$ as $(p,q)$-forms
-- clearly $H^\bullet(X, \CC)$ possesses a basis consisting of such forms.  The
cup/wedge product on cohomology furnishes this vector space with the
structure of a bigraded $\CC$-algebra.  Finally, integration of forms
induces a trace on the algebra; in terms of a basis element $\omega$, 
\[
\text{tr}(\omega) = 
\begin{cases}
  \displaystyle\int_X \omega & \omega \in H^n(X, \ext n T_X^\vee)\\
  0 & \text{otherwise.}
\end{cases}
\]
The pairing $(\alpha, \beta) \mapsto \text{tr}(\alpha\wedge\beta)$ induced by this trace
is a non-degenerate bilinear form satisfying \((\alpha, \beta \wedge\gamma)
= (\alpha \wedge\beta, \gamma),\)
so that $H^\bullet(X, \CC)$ is a bigraded Frobenius algebra.

\subsection*{Physics}

The relationship between ordinary cohomology and quantum cohomology may be
elucidated by appealing to physics -- in particular to a topologically-twisted
(2,2) nonlinear sigma model of maps $\PP^1 \to X$.  
Of the many intriguing aspects of this quantum field theory, we will be
most interested in its algebra of massless supersymmetric
operators\footnotemark.  Using
elementary physics arguments, one identifies a basis for the set of
such operators that may be set into one-to-one correspondence with $(p,q)$-forms on
$X$. 

\footnotetext{More precisely, it is the algebra of local, scalar,
supersymmetric operators \cite{Witten:1991zz}.}

The (2,2) supersymmetry of the model forces the product of two massless
supersymmetric operators to be massless and supersymmetric.  The particular
form of the product obtains by considering three-point correlation
functions in the quantum field theory: the quantum product of two massless
operators $\sheaf O_1$ and $\sheaf O_2$ is defined to be the unique
operator $(\sheaf O_1 * \sheaf O_2)$ such that for all massless operators
$\sheaf O_3$,  
\begin{equation}
  \left\langle 
  \mathds 1 (\sheaf O_1 * \sheaf O_2) \sheaf O_3 \right\rangle  = 
  \left\langle \trip \right\rangle.
  \label{eq:quantumProduct}
\end{equation}
Here, $\mathds 1$ denotes the
operator corresponding to $1 \in H^0(X, \ext 0 T_X^\vee)$.  Such a
correlation function is computed using the instanton expansion  
\begin{equation}
  \left\langle \trip \right\rangle = \sum_{\beta
  \in H_2(X,\ZZ)} \left\langle \trip 
  \right\rangle_\beta q^\beta.
  \label{eq:instanton}
\end{equation}
Although they have intrinsic meaning in physics, we will consider the
expressions $q^\beta$ to comprise a set of formal variables endowed with
the structure of a monoid via the product $q^\alpha q^\beta =
q^{\alpha+\beta}$.  We denote by $\CC\llbracket \underline q \rrbracket$
the ring of formal power series with complex coefficients in these
variables -- one sometimes insists on convergence or other properties, but
such subtleties are beyond the scope of this review.

Mathematically, one defines the expression $\left\langle \trip
\right\rangle_\beta$ as the Gromov-Witten invariant\footnotemark~$\langle
I_{0,3,\beta}\rangle(\omega_1, \omega_2,\omega_3)$, where $\omega_i$ are
the forms corresponding to the operator $\sheaf O_i$. 
\footnotetext{See equation 7.4 of \cite{Cox:2000vi} for a precise
definition of Gromov-Witten invariants.}
Physically, one says that $\left\langle \trip \right\rangle_\beta q^\beta$
denotes the contribution of instantons of degree $\beta$ to the correlation
function $\left\langle \trip \right\rangle$.  This expression is morally
the integral of induced forms on some compactification $\modcpt$ of the
moduli space of holomorphic maps $f:\PP^1\to X$ of class $\beta =
f_*[\PP^1]$.  We will write the induced forms schematically using maps 
\begin{equation}
  \zeta_\beta \takes {H^p(X, \ext q T_X^\vee)}{H^p\left(\modcpt, \ext q
  T_\modcpt^\vee\right)}.
  \label{eq:inducedMap}
\end{equation}
If $\omega_i$ are the forms corresponding to operators $\sheaf O_i$, modulo
the subtleties of obstruction bundles we have that 
\[
\left\langle \trip \right\rangle_\beta = 
\int_\modcpt \zeta_\beta(\omega_1) \wedge \zeta_\beta(\omega_2) \wedge
\zeta_\beta(\omega_3).
\]
Depending on the compactification, there may be more than one such map --
in the case of the stable maps compactification, pullbacks via evaluation
maps play the r\^ole of $\zeta_\beta$.  For toric varieties, one often uses
the Morrison-Plesser compactification \cite{Morrison:1994fr} wherein -- as
indicated in \eqref{eq:inducedMap} -- one map suffices for each $\beta$.

The three-point correlation functions in \eqref{eq:instanton} induce a
non-degenerate bilinear pairing $(\omega_1, \omega_2) = \langle \mathds 1
\sheaf O_1 \sheaf O_2 \rangle$ on the unital algebra $\bigoplus_{p,q}
H^p(X, \ext q T_X^\vee)\llbracket \underline{q} \rrbracket$, leading to the
following definition.
\begin{definition}
  The \emph{quantum cohomology} of $X$ is the Frobenius algebra
  \[
  QH^\bullet(X) := \bigoplus_{p,q} H^p(X, \ext q T_X^\vee)\llbracket \underline{q} \rrbracket,
  \]
  with the product and bilinear pairing induced by (2,2) three-point 
  functions.
\end{definition}
Here, the (2,2) correlation functions are defined either via
Gromov-Witten invariants or as correlation functions in the quantum field
theory, depending on whether your tastes tend to the mathematical or
to the physical. 

\begin{example}
The ``classical sector'' is the set of maps homotopic to a
point, $\beta = 0$, and the moduli space of such maps is simply $X$ itself.
Thus, in this sector, the quantum product reduces to the wedge product on forms; ordinary
cohomology is the ``classical limit'' of quantum cohomology.  This sector
may be isolated 
by setting $\underline q = 0$. For example,
  the ordinary and quantum cohomology of $\PP^n$ are respectively
  \[
  \begin{split}
	 H^\bullet(\PP^n, \CC) &\cong \frac{\CC[H]}{\left\langle H^{n+1} \right\rangle},\\
	 QH^\bullet(\PP^n) &\cong \frac{\CC[H]\llbracket q\rrbracket}{\left\langle H^{n+1} -q \right\rangle }.
  \end{split}
  \]
  Here $H$ denotes the hyperplane class. For $\PP^n \times \PP^m$, the
  equivalent expressions are 
  \begin{equation}
	 \begin{split}
		H^\bullet(\PP^n\times \PP^m, \CC) &\cong \frac{\CC[H_1,H_2]}{\left\langle H_1^{n+1},
		H_2^{m+1} \right\rangle},\\
		QH^\bullet(\PP^n\times \PP^m) &\cong \frac{\CC[H_1,H_2]\llbracket
		q_1,q_2\rrbracket}{\left\langle H_1^{n+1}
		-q_1, H_2^{m+1} -q_2 \right\rangle }.
	 \end{split}
	 \label{eq:qcohomPnPm}
  \end{equation}
\end{example}

\section*{Quantum Sheaf Cohomology}

As in our study of the passage from ordinary cohomology to quantum
cohomology, we first consider the ``ordinary sheaf cohomology'' -- in
particular that of an omalous bundle $\sheaf E \to X$.  Here, by ordinary
sheaf cohomology we mean 
cohomology valued in polysections,
\begin{equation}
  \bigoplus_{p,q} H^p(X, \ext q \sheaf E^\vee).
  \label{eq:classicalsheaf}
\end{equation}
Again by a slight abuse of language, we
will refer to elements of $H^p(X, \ext q \sheaf E^\vee)$ as
$(p,q)$-forms -- clearly the vector space \eqref{eq:classicalsheaf} possesses a basis
consisting of such forms, and the cup/wedge product furnishes it with the
structure of a bigraded $\CC$-algebra.  

The trace on this algebra is slightly more subtle and follows from the omality
of $\sheaf E$.  In particular, one uses the existence of an isomorphism
$\psi\takes{H^n(X, \ext {k}\sheaf E)}{H^n(X, \omega_X)}$ to define, for
a basis element $\omega$, 
\begin{equation}
  \text{tr}(\omega) = 
  \begin{cases}
	 \displaystyle\int_X \psi(\omega) & \omega \in H^n(X, \ext k \sheaf E^\vee)\\
	 0 & \text{otherwise.}
  \end{cases}
  \label{eq:sheafTrace}
\end{equation}
The pairing $(\alpha, \beta) \mapsto \text{tr}(\alpha\wedge\beta)$ induced by this trace
endows $\bigoplus_{p,q} H^p(X, \ext q \sheaf E^\vee)$ with the structure of a bigraded Frobenius algebra.

\subsection*{Physics}

To understand the relationship between sheaf cohomology and quantum sheaf
cohomology we again appeal to physics -- in particular a
topologically-twisted (0,2) nonlinear sigma model of maps $\PP^1 \to X$.  A
recent physics review of this and related models may be found in
\cite{McOrist:2010ae}.  We will again be most interested in its algebra of
massless supersymmetric operators\footnotemark.  The same elementary
physics arguments used for the (2,2) theory identify a basis for this set
of operators that may be placed into one-to-one correspondence with
$(p,q)$-forms (that is, elements of \eqref{eq:classicalsheaf}), and the
quantum product of two massless operators is defined using three-point
correlation functions of the (0,2) in analogy to \eqref{eq:quantumProduct}.
Unlike the (2,2) case, however, there is no mathematical definition of
$\left\langle \trip \right\rangle_\beta$ in a (0,2) theory so the following
definition is purely physical.

\footnotetext{As explained in \cite{Adams:2005tc}, we are actually interested in
local, scalar, supersymmetric operators with vanishing holomorphic
conformal weight, but for continuity we will refer to them as 
massless supersymmetric or simply massless.}
\begin{definition}
  The \emph{quantum sheaf cohomology} of an omalous bundle $\sheaf E \to X$ 
  is the Frobenius algebra
  \[
  QH^\bullet(X, \sheaf E) := \bigoplus_{p,q} H^p(X, \ext q \sheaf E^\vee)\otimes \CC\llbracket \underline q \rrbracket
  \]
  with the product and bilinear pairing induced by (0,2) three-point
  functions.
\end{definition}
As in the case of ordinary quantum cohomology, the classical limit of
quantum sheaf cohomology is precisely the ordinary sheaf cohomology with
the Frobenius structure induced by \eqref{eq:sheafTrace}.
Unlike the case in (2,2) theories, (0,2) supersymmetry is \emph{not} enough
to guarantee that the product of massless operators is massless: one needs
to work harder to show that the algebra closes in the set of all
operators.

\subsection*{Existence} 

The (modern) history of quantum sheaf cohomology begins with the
observation in \cite{Adams:2003zy} of an analogue of $QH^\bullet(X)$ for $(0,2)$
theories.  Therein, the quantum sheaf cohomology of a one-parameter family
of deformations of the tangent bundle of $\PP^1 \times \PP^1$ was computed
using a conjectured form of mirror symmetry for $(0,2)$ models.  Their
calculations were confirmed in a sheaf-cohomology-based computation by Katz
and Sharpe \cite{Katz:2004nn}.   Inspired by these results, Adams, Distler,
and Ernebjerg \cite{Adams:2005tc} gave a physics definition of quantum
sheaf cohomology and found a physics proof of two sufficient conditions for
its existence.  We restate these conditions here as conjectures.

\begin{conjecture}
  Let $\sheaf E$ and $\sheaf E'$ be omalous elements of a family of bundles
  $U$.  Let $\gamma\!:\![0,1] \to U$ continuous with
$\gamma(0) = \sheaf E$,
$\gamma(1) =\sheaf E'$,
$\gamma(t)$ omalous for all $t \in [0,1]$. 
  Then $QH^\bullet(X, \sheaf E)$ exists iff $QH^\bullet(X, \sheaf E')$ exists.
\end{conjecture}

\begin{conjecture}
  If $\sheaf E \to X$ is omalous and rk $\sheaf E < 8$, then $QH^\bullet(X, \sheaf
  E)$ exists.
\end{conjecture}

Since $QH^\bullet(X, T_X) = QH^\bullet(X)$, 
the former condition implies the existence of quantum sheaf cohomology for
all omalous one-parameter families of tangent-bundle deformations.
The latter is likely an artefact of the technique used in the physics proof
-- there are no known examples of omalous bundles of rank eight or
higher for which the massless operators fail to close under the quantum
product, and there are no physical reasons to expect such a bundle to
exist.

\subsection*{Computation example}
Although there is no definition for the invariants $\left\langle \trip
\right\rangle_\beta$, a number of physics-inspired techniques exist to
compute them 
when the omalous bundle is a deformation of the tangent bundle
of a toric variety
\cite{Katz:2004nn, Guffin:2007mp, McOrist:2007kp} or a complete
intersection therein\cite{McOrist:2008ji}.
One of the advantages of using a toric variety $X$ is that deformations of
$T_X$ are easily obtained by deforming the Euler exact sequence:
\[
\begin{tikzpicture}[description/.style={fill=white,inner sep=2pt},baseline=(current bounding box.center)]
  \matrix (m) [matrix of math nodes, row sep=1em, column sep=1.5em, text height=1.5ex, text depth=0.25ex]
  {
  0 & \mathcal O_X^r & \displaystyle\bigoplus_{\rho \in \Delta}
  \mathcal O_X(D_\rho) & T_X & 0.\\
  };
  \draw[->] (m-1-1) -- (m-1-2);
  \draw[->] (m-1-2) -- (m-1-3) node [midway,above] {\small $E_0$};
  \draw[->] (m-1-3) -- (m-1-4);
  \draw[->] (m-1-4) -- (m-1-5);
\end{tikzpicture}
\]

Here, $r$ is the rank of the Picard group, $\Delta$ denotes the set
of torus-invariant divisors corresponding to one-cones in the fan of $X$, and $E_0$ is a collection of sections of $\sheaf
O_X(D_\rho)$, which are toric analogues of $\sheaf O_{\PP^n}(1)$.  Taking
$X = \PP^1 \times \PP^1$, for example, the sequence becomes
\[
\begin{tikzpicture}[description/.style={fill=white,inner sep=2pt},baseline=(current bounding box.center)]
  \matrix (m) [matrix of math nodes, row sep=1em, column sep=1.5em, text height=1.5ex, text depth=0.25ex]
  {
  0 & \mathcal O_X^2 & \mathcal O_X(1,0)^2 \oplus \mathcal O_X(0,1)^2 &
  T_X & 0, 
  \\
  };
  \draw[->] (m-1-1) -- (m-1-2);
  \draw[->] (m-1-2) -- (m-1-3) node [midway,above] {\small $E_0$};
  \draw[->] (m-1-3) -- (m-1-4);
  \draw[->] (m-1-4) -- (m-1-5);
\end{tikzpicture}
\]
where the map is 
\[
E_0 = \left( \begin{matrix}
  x_0 & 0 \\
  x_1 & 0 \\
  0 & x_2 \\
  0 & x_3 
\end{matrix}\right).
\]

A deformation of $T_X$ may be obtained by choosing a different collection
of sections for the map.  For example, selecting the map to be 
\begin{equation}
E = \left( \begin{matrix}
  x_0 & \epsilon_1 x_0 + \epsilon_2 x_1 \\
  x_1 & \epsilon_3 x_0 \\
  \gamma_1 x_2 + \gamma_2 x_3 & x_2 \\
  \gamma_3 x_2 & x_3 
\end{matrix}\right)
  \label{eq:deformation}
\end{equation}
as in \cite{Guffin:2007mp} gives a convenient basis for the space of
deformations of the tangent bundle ($\epsilon_i, \gamma_i \in \CC$).
Therein, several of the invariants $\left\langle \trip \right\rangle_\beta$
were computed for the bundle $\sheaf E \to \PP^1 \times \PP^1$ defined as
the 
cokernel of \eqref{eq:deformation}.  These were then used to deduce the
quantum sheaf cohomology of $\sheaf E$;
\begin{equation}
  QH^\bullet(\PP^1\times\PP^1, \sheaf E) \cong \frac{\CC[\psi,\wt\psi]\llbracket
  q_1,q_2\rrbracket}{\left\langle
  \begin{matrix}
	 \psi^2 + \epsilon_1 \psi \wt\psi - \epsilon_2 \epsilon_3 \wt\psi^2 -q_1, \\
	 \wt\psi^2 + \gamma_1 \psi \wt\psi - \gamma_2 \gamma_3 \psi^2
	 -q_2\phantom{,,} 
  \end{matrix}
  \right\rangle }.
  \label{eq:qscohomP1P1}
\end{equation}
These computations were confirmed in \cite{McOrist:2007kp} using physics
techinques.  
Note that as $\epsilon_i, \gamma_i \to 0$, the quantum sheaf
cohomology in \eqref{eq:qscohomP1P1} limits to the ordinary quantum
cohomology of $\PP^1 \times \PP^1$ in \eqref{eq:qcohomPnPm}.

This material is based on work supported by the National Science Foundation
under DMS grant no.~0636606.
\bibliography{/work/all}
\end{document}